\documentclass{tran-l}
\usepackage{tipa}
\usepackage{amssymb}
\usepackage{amsmath}
\usepackage{tikz}
\usepackage{mathtools}
\usepackage{tikz}
\usepackage{tikz-cd}
\usetikzlibrary{matrix}
\usepackage{float}
\usepackage{graphicx}
\usepackage{subcaption}
\usepackage{xcolor}
\usepackage{enumitem}
\usepackage[utf8]{inputenc}

\newtheorem{thm}{Theorem}[section] % 1st argument is your name for it
\newtheorem{lem}[thm]{Lemma}     % 2nd argument is what is printed
\newtheorem{cor}[thm]{Corollary}
\newtheorem{prop}[thm]{Proposition}
\theoremstyle{definition}
\newtheorem{defn}[thm]{Definition}
\newtheorem{rem}[thm]{Remark}
\newtheorem{exa}[thm]{Example}
\numberwithin{equation}{subsection}
\newcommand{\eps}{\varepsilon}

\newcommand{\de}{\delta}

\newcommand{\si}{\sigma}
\newcommand{\SI}{\Sigma}

\newcommand{\II}{\mathbb{I}_M}

\newcommand{\R}{\mathbb{R}}

\newcommand{\FM}{\mathfrak{M}}

\newcommand{\ds}{\mathfrak{ds}}

\begin{document}

 \title[]
 {Metric Spaces and Sparse Graphs}

\author{  J.M. Alonso }

\address{IMASL-CONICET. Universidad Nacional de San Luis. Ej\'ercito de Los Andes 950. 5700 San Luis. Argentina}

\email{jmalonso@unsl.edu.ar, jm31415ac@gmail.com}

%\thanks{}


%\subjclass{92E10 (primary), 68R10}


% -----------------------------------------------------------

\begin{abstract}
Many concrete problems are formulated in terms of a finite set of points in $\R^N$ which, via the ambient Euclidean metric, becomes a finite metric space. To obtain information from such a space, it is often useful to associate a graph to it, and do mathematics on the graph, rather than on the space. 

Connected graphs become finite metric spaces (rather, their set of vertices) via the path-metric. We consider different types of connected graphs that can be associated to a metric space. In turn, the metric spaces obtained from these graphs recover, in some cases, identically the initial space, while they, more often, only ``approximate'' it. In this last category, we construct a \emph{connected sparse} graph, denoted $CS$, that seems to be new. We show that the "general case" for $CS$ is to be a tree, a result of clear practical importance.
\end{abstract}

% -----------------------------------------------------------
\maketitle
% -----------------------------------------------------------

\section*{Introduction} 
In \cite{aloet21} we presented a novel mathematical	approach to study the 
nature of the glass and the glass transition. The method used can be summarized as follows. From simulations at different temperatures, we obtained the cartesian coordinates of about 2,200 silicon and oxygen atoms. For seven chosen temperatures we took, at very short time intervals, 101 ``snapshots''. The idea was to obtain information from these 707 sets of points in $\R^3$ by computing their \emph{finite dimension} \cite{alo16, alo15}, a version, adapted to finite sets, of the classical Hausdorff dimension. Finite dimension, denoted $\dim_f$, is defined (only) on finite metric spaces and its values are highly non-trivial, in contrast to the classical Hausdorf dimension. Unfortunately, it is not effectively computable on such general and large sets. This paper studies in detail the solution to this problem implemented in \cite{aloet21}.

In the old days when homology groups where defined only for simplicial complexes, Vietoris \cite{vie} associated a simplicial complex to a metric space, the Vietoris complex, thus extending the definition of homology groups to arbitrary metric spaces. This construction was rediscovered by Rips in the 1980's in the context of Gromov's theory of discrete hyperbolic groups \cite{gro1} (under the name Rips complex) where it was used for purposes other than computing homology groups (see e.g. \cite{alo+}, \cite{alo}). Nowadays it is usually called the Vietoris-Rips complex. 

The Vietoris-Rips complex should be considered an ``approximation'' of the metric space. In some instances this can be made precise: if the metric space is a closed Riemannian manifold and the Vietoris-Rips complex corresponds to a sufficiently small cut-value, then the manifold and the complex are homotopy equivalent \cite{hau, lat}.

We recycled this idea by considering only the 1-skeleton of the Vietoris-Rips complex, a graph, and by taking advantage of the finiteness of our spaces. The result is the \emph{minimum connected graph}, $MCG$, defined in Section \ref{s:mcg}. Thus, going back to the glass simulations, we reduced the problem from computing the finite dimension of the ``snapshots'' to computing the finite dimension of the ``approximating''  $MCG$s. It turns out that $MCG$ is not completely satisfactory for effective computation because, while it is connected (more on this below), it usually contains many ``superfluous'' edges, superfluous in the sense that they are not necessary to assure connectedness. Searching for a streamlined version of $MCG$, we came to define the \emph{connected sparse graph}, $CS$, associated to a finite metric space. This graph is the main object of study of this paper. It turns out that, for the 707 spaces we studied in \cite{aloet21}, the $CS$s are trees\footnote{This fact seemed too good to be a mere coincidence (cf. Section \ref{s:gen_case}).} while the $MCG$s contain about 3 times as many edges as vertices. Thus, we ``approximated'' the glass simulations by the $CS$s and could then effectively compute their finite dimension, see the results in \cite{aloet21}.

To extend the applicability of $CS$ to other situations, we start with a finite metric space $(M,d)$ (usually $M\subset \R^N$ and $d$ is the restriction of the Euclidean metric) and define a graph $G(M)$ associated to it, in such a way that its set of vertices equals $M$. So $G(M)$ is a graph ``overlaid'' on top of $M$. We further require that $G(M)$ be connected so that, from $G(M)$ we can obtain a new metric structure on $M$, by using the path-metric (see Section \ref{s:grtom}). Sometimes, this procedure reproduces the original $(M,d)$, more often, it only ``approximates'' it. We insist in requiring $G(M)$ to be connected because, by definition, $\dim_f(G(M))$ is the finite dimension of $M$ (the set of vertices) with the path-metric obtained from $G(M)$. The graphs $CS$ and $MC$ are only two of many possible choices for $G(M)$. As for notation, since we only consider graphs, we drop the $G$ and simply write $CS$ and $MC$. If we need to stress the space $CS$ [or $MC$] comes from, we write $CS(M)$ [or $MC(M)$], Schematically, the plan can be pictured as follows:
\[
(M,d) \quad   \xmapsto[]{\qquad} \quad G=G(M)\quad \xmapsto{\qquad} \quad \textrm{compute on } G: 
   \begin{dcases}
     \textrm{Finite Dimension}\\
     \textrm{Harmonic Analysis} \\
     \textrm{Persistent Homology}\\
     \cdots \cdots \cdots \cdots  \cdots \cdots
   \end{dcases}
\]
ie. extract information from $(M,d)$ by applying your favourite tool to $G(M)$.

The construction of $CS(M)$ is unstable: a slight perturbation of the vertices, ie. of the points of $M$, can produce very different types of graphs (cf. Example \ref{exa:unst}); in particular, graphs that are far from being sparse. To study this phenomenon we define a new metric space, $(\FM_m, d_B)$, where $\FM_m$ is the set of all finite sets $M\subset \R^N$ of cardinality $m$, with the induced Euclidean distance. Thus, $\FM_m$ is the space of all metric subspaces of $\R^N$ of fixed cardinality $m$. We use $\FM_m$ to show that, if $M$ is contained in some Euclidean space and we are willing to perturb it slightly, then we can make $CS(M)$ a tree, a fact of clear practical importance, since a slight perturbation of $M$ is within empirical error. Thus, for $M\subset \R^N$, the  ``general case'' is for $CS(M)$ to be a tree, thus justifying the name.

The contents of the paper are as follows. In Section \ref{s:mcg} we recall several well-known facts about graphs, establish notation, and define the minimum connected graph $MC(M)$. The connected sparse graph $CS(M)$ is defined in Section \ref{s:csg}. In Section \ref{s:space_of_sp} we define $\FM_m$ and its subspace $\ds_m$ of \emph{distance separated spaces}, and prove that  $\ds_m$ is open and dense in $\FM_m$. Section \ref{s:gen_case} contains the main result of the paper, that the ``general case" is for $CS(M)$ to be a tree. Both $MC(M)$ and $CS(M)$  are invariant under isometries of $M$, as  shown in Section \ref{s:inv}.  Graphs $G(M)$ that exactly reproduce the original $M$ are discussed in Section \ref{s:intri}, and the relation among the different $G(M)$ we have considered here, is discussed in Section 7.

\section{The Minimum Connected Graph}\label{s:mcg}

The path-length defined on connected graphs is a well-known method, recalled in Section \ref{s:grtom} below, to give the set of vertices a structure of metric space. The metric spaces thus obtained (called \emph{intrinsic-I} in Section \ref{s:intri} below) depend only on the isomorphism class of the graph \cite{alo16}. This method can be generalized to graphs with positive weights, not necessarily equal for all edges. The corresponding spaces (called \emph{intrinsic-II} in Section \ref{s:intri}) are, however, a medley between combinatorics and the reals, not uniquely determined by the combinatorics. Beyond this, finite metric spaces are \emph{extrinsic} (see Section \ref{s:intri}). This classification of finite metric spaces is particularly interesting from the point of view of finite dimension. In short, $\dim_f$ is effectively computable for intrinsic-I spaces, one can sometimes use the graph structure in the case of intrinsic-II spaces, and it is usually not effectively computable for extrinsic spaces. 

Finite dimension has been applied in situations where there is an underlying graph with constant weights (glycans \cite{aloet18}) and in other situations where there is a natural graph, but the weights are not constant (plants and roots \cite{aloet19}). In the case of glasses \cite{aloet21} it is necessary to mediate $CS(M)$ and compute its finite dimension in order to obtain information about $M$.

\subsection{From graphs to metric spaces} \label{s:grtom}

We now recall the construction of a metric space from a graph and establish notation. In this paper, graph means an undirected, finite, connected, simple graph without loops; equivalently, the 1-skeleton of a finite, connected, simplicial complex. Explicitly, a graph $G=(V,E)$ consists of a set $V$ of vertices or nodes, and a set $E$ of edges, where edges are of the form $\{ v,v'\}$, for $v\neq v'$ in $V$; when $\{ v,v'\}\in E$, we say that $v$ and $v'$ are \emph{adjacent}. A \emph{path} $p(v,v')$ connecting $v$ to $v'$, is a sequence $v_0,\dots,v_n$ of distinct nodes such that $v_i$ and $v_{i+1}$ are adjacent and $v_0=v, v_n=v'$. We call G \emph{connected} when any two  nodes are connected by a path. We also write $p(v,v')=e_1e_2\cdots e_n$, where the $e_i\in E$ are of the form $e_1=\{v_0,v_1\}, e_2=\{v_1,v_2\},\cdots, e_n=\{v_{n-1},v_n\}$. For such $p$, we define the \emph{count} of $p$, $c(p)$, to be $n$. A \emph{cycle} in $G$ is a sequence $v_1,\dots,v_n,\, n\geq 3,$ of distinct nodes such that $v_i$ and $v_{i+1}$ are adjacent (for $i=1,\dots,n-1$) and $v_1$ and $v_n$ are adjacent.

Edges $e\in E$ have weights $w_{e}>0$, which we interpret as lengths. A path $p(v,v')$ as above has \emph{length} $\sum_{i=1}^k w_{e_i}$. From $G$ we derive a metric space $(V,d_G)$, where $d_G(v,v')$ is the minimum length of all paths connecting $v,v'$. A path $p(v,v')$ is a \emph{geodesic} if its length equals $d_G(u,v)$. We abuse language and call $G$ a metric space, when in fact we mean $(V,d_G)$. In this way, graphs become a rich source of examples of intrinsic (to be defined later on) finite metric spaces.

\subsection{From metric spaces to graphs} \label{s:mtogr}
In this paper we go in the opposite direction, from metric spaces to graphs, from $M$ to $G(M)$. It is natural to give $G(M)$ one of the following two types of weights: all constant, equal to 1, or $w_e=d(x,y)$ for every edge $e=\{x,y\}$. In the former case we use the notation $G(M)$, in the latter, $G(M)^w$ and, for the corresponding path-metric, $d_{G(M)}$ and $d_{G(M)}^w$. Thus, it is only in the second case that we make the weights explicit.

%\section{The Minimum Connected Graph}\label{s:mcg}
We now construct, given a finite metric space $(M,d)$, the graph $MC(M)$. For this purpose, instead of the whole Vietoris-Rips complex, we consider its 1-skeleton, which is a simple graph, and then focus on the smallest such that is connected. Here are the details. 

Consider the \emph{set of distances} $DM:=\{d(x,y)\in\R | x,y\in M\}=\{r_0,r_1,\dots,r_m\}$, where $0=r_0<r_1<\cdots<r_m$, and $r_m$ is the diameter of $M$. For $r\in DM$, let $G_t$ denote the graph with vertices $V_r=M$, and edges $E_r:=\{\{x,y\}\in M^2 | d(x,y)\leq r\}$. We say that $r$ is the \emph{cut value} that defines $G_r$. Observe that $G_{r_{i+1}}$, is obtained from $G_{r_i}$ by adding precisely those edges $\{x,y\}$ for which $ d(x,y) = r_{i+1}$. 

Note that $G_{r_0}$ is a graph with no edges and $|M|$ connected components, one for each point of $M$. At the other extreme, $G_{r_m}=K_{|M|}$, the complete graph on $|M|$ vertices. Somewhere between these two extremes we find $MC=MC(M)$, the minimum connected graph of $M$. More formally, there is a smallest $i, 0<i\leq m$, such that $G_{r_i}$ is connected, and we set $MC(M):=G_{r_i}$. Of course, we can also use the distances in $M$ as weights, and obtain $MC(M)^w$.

\section{The Connected Sparse Graph} \label{s:csg}
The $MC$s contain many ``superfluous'' edges (particularly if we are interested in obtaining a connected graph) resulting, in general, in a large graph. The definition of the connected sparse graph, $CS$, aims at obtaining a sparse graph by including only  those edges that are ``necessary'' to obtain a connected graph. Since we also want $CS$ to be unique, we make no choices and, in doing so, may be forced to add a few superfluous edges (superfluous in the sense that adding fewer edges one might still obtain a connected graph). The definition is inductive: starting with $M$, new edges are added at each step, obtaining an increasing sequence of graphs:
\begin{displaymath}
S_0 \subset S_1 \subset \cdots 
\end{displaymath}

\noindent
STEP 0. $S_0$ is a graph with no edges, $V(S_0)=M$, and $E(S_0)=\emptyset$. The set of components of $S_0$ is $K_0
=\{\{x\}| x \in M\}$, which we identify with $M$ by setting $\{x\}=x$. Thus $S_0$ has $|K_0|=|M|$ components.

\noindent
STEP 1. The function $\nu_0:M\to \R_+$ (where $\R_+$ denotes the positive reals), defined by 
\begin{displaymath}
\nu_0(x):= \min_y \Big\{d(x,y)| y\in M, y\neq x \Big\}, 
\end{displaymath}
gives the distance to the nearest neighbours of a point. Note that we can consider $\nu_0$ as a function defined on $K_0$, by setting $\nu_0(\{x\}):= \nu_0(x)$. Let

\begin{displaymath}
E(S_1):= \bigcup_{x\in M} \Big\{ \{x,y\}| d(x,y)= \nu_0(x) \Big\}.
\end{displaymath}
That is, for every vertex $x$, we add an edge connecting $x$ to each one of its nearest neighbours.
Clearly, $S_0\subset S_1$.

The graph $S_1$ is not connected in general. Let $K_1$ denote the set of its components. If $|K_1|=1$, we stop here and define $CS(M):= S_1$. If $|K_1|>1$, we proceed to Step 2. Note that $|K_1|\le |M|/2$, since every $x \in M$ has a nearest neighbour.

\noindent
STEP $i+1$. In this inductive step, we assume the graphs $S_0\subset S_1 \subset \cdots\subset S_{i}$ are constructed, that   $K_{i}$ is the set of connected components of $S_{i}$, and that $1<|K_{i}| \le |M|/2^{i}$.

Let $k,k'\in K_{i}$ be distinct components. The function $\nu_{i}: K_{i} \to  \R_+$ , defined by 
\begin{displaymath}
\nu_{i}(k):= \min_{x,y} \Big\{d(x,y)| x\in k, y\in k', k\neq k'\Big\}=  \min_{x,y} \Big\{d(x,y)| x\in k, y\notin k\Big\},
\end{displaymath}
gives the distance from $k$ to its nearest components. For $x\in M$, we let $[x]_i$ denote the component of $x$ in $S_i$. Let
\begin{displaymath}
F(S_{i+1}):= \bigcup_{k\in K_{i}} \Big\{ \{x,y\}| d(x,y)= \nu_{i}(k), x\in k, y\notin k \Big\}.
\end{displaymath}
and define $E(S_{i+1}):=E(S_{i}) \cup F(S_{i+1})$. Clearly, $S_{i}\subset S_{i+1}$.

Let $K_{i+1}$ denote the set of components of $S_{i+1}$. If $|K_{i+1}|=1$, we stop here and define $CS(M):= S_{i+1}$. If $|K_{i+1}|>1$, we proceed to Step ${i+2}$. Note that $|K_{i+1}|\le |K_{i}|/2$, since every $k \in K_{i}$ has a nearest connected component. Hence, $|K_{i+1}|\le |M|/2^{i+1}$.

Since $|M|$ is finite, after a finite number of steps we obtain a connected graph, and this we define to be $CS(M)$. Also, note that no choices are made at any step, so that $CS(M)$ is uniquely defined. 

\begin{rem} Since $CS(M)$ depends completely on the distance of $M$, it  reflects properties of the metric space, a fact we express by saying that $CS(M)$ approximates $(M,d)$. Of course, $CS(M)^w$ is an even better approximation to $(M,d)$.
\end{rem}

\begin{exa} \label{exa:unst}
For $x\leq 0$, consider the set $M_x$ of the following $n+1$ points of $\R^n$, $a_x=x(1,1,\dots,1)$, and $e_i=(0,\dots,0,1,0,\dots,0),\,(1\leq i\leq n)$, where 1 is the $i$th-coordinate. With the Euclidean distance, $d(e_i,e_j)=\sqrt{2}$, for $i\neq j$, and: 

\[ 
d(a_x,e_i) =
\begin{cases}
 1    & \text{for } x=0, \\
 \sqrt{2}     & \text{for }x=(1-\sqrt{1+n})/n.
\end{cases}
\]
Moreover, $d(a_x,e_i)>\sqrt{2}$ for $x<(1-\sqrt{1+n})/n$, and $1 <d(a_x,e_i)< \sqrt{2}$, for $(1-\sqrt{1+n})/n< x <0$. It follows that $CS(M_x)=St_{n+1}$, the \emph{star graph} on $n+1$ vertices, for $(1-\sqrt{1+n})/n < x \leq 0$, and $CS(M_x)=K_{|M_x|}$, the complete graph on $|M_x|$ vertices, for $x\leq (1-\sqrt{1+n})/n$. Finally, note that $CS(M_x)=MC(M_x)$.
\end{exa}

The following technical proposition will be needed later.

\begin{prop} \label{p:tech}
Let $x,y\in M$, and let $[x]_i\in K_i$ denote the unique component (in $S_i=S_i(M)$) to which $x$ belongs. Then 
\begin{enumerate}
  \item $E(S_{i+1})=\bigcup_{j=1}^{i+1} F(S_{j})$, where the union is disjoint; in other words, the family $\{F(S_j)\}_{j=1}^{i+1}$ is a partition of $E(S_{i+1})$.
  \item $\{x,y\} \in F(S_{i+1})$ if and only if
\begin{displaymath}
\Big([x]_i\neq [y]_i \Big) \textrm{ and }  d(x,y) = \max \Big\{ \nu_i([x]_i), \nu_i([y]_i) \Big\}.
\end{displaymath}
\end{enumerate}

\end{prop}

\begin{proof} Both conditions follow directly from the inductive definition of the $S_j$. 
\end{proof}

%XXXXXXXXXXXXXXXXXXXXXX
%XXXXXXXXXXXXXXXXXXXXXX
%XXXXXXXXXXXXXXXXXXXXXX

\section{The space of finite metric spaces of fixed cardinality.} 
\label{s:space_of_sp}
In this section and the next, all metric spaces are subsets of Euclidean space, metrised as subspaces of $\R^N$. Let $\FM_m$ denote the set of subsets $M\subset \R^N$ of cardinality $m$, and let $d$ denote the Euclidean distance in $\R^N$.
We define a distance $d_B$ in $\FM_m$:

\begin{defn}
For $M,M'\in \FM_m$ and a bijection $f:M\to M'$, define $\textrm{sep }f$, the \emph{separation} of $f$, by:
\[
\textrm{sep }f := \max \{d(x,f(x))| x\in M\},
\]
and a distance
\[
d_B(M,M'):=\min\{\textrm{sep }f | f\in \textrm{B}(M,M')\},
\]
where $\textrm{B}(-,-)$ denotes the set of bijections.
\end{defn}

It is easy to see that $\textrm{sep }f = \textrm{sep }f^{-1}$ and $\textrm{sep }gf\leq \textrm{sep }g + \textrm{sep }f$, for composable bijections $g,f$. This can be used to show that $d_B$ is a metric.

\begin{defn}
Recall from Section \ref{s:mcg} the set $DM:= \{ d(x,y)| x,y\in M\}$ of distances of $M$, and let $DM_+:=DM\setminus \{0\}$, denote the set of \emph{positive distances}. The \emph{mesh} of $DM$, $\de=\de(DM)>0$ is, by definition, the smallest positive distance between the points of $DM$, i.e. $\min\{|r-r'||r,r'\in DM, r\neq r'\}$. Note that $\de(DM)\leq \de(M)$.
\end{defn}

\begin{defn}
For $x_i\in M, r\in DM$, let $\SI(x_i,r):=\{y \in M | d(x_i,y)=r \}$ denote the sphere centered at $x_i$, with radius $r$. By the choice of $x_i$ and $r$, $\SI(x_i,r)\neq \emptyset$ .
\end{defn}

\begin{rem}
We reserve the notation $S(x,r):=\{y\in \R^N| d(x,y)=r\}$ for spheres centred at arbitrary points of $\R^N$, with any positive real radius. Obviously, $M\cap S(x,r)$ can be empty, even if $x\in M$.
\end{rem}

\begin{defn}
$(M,d)$ is \emph{distance separated} if $DM_+$ has exactly $m(m-1)/2$ elements, i.e. $DM_+$ has the largest possible cardinality. Equivalently, all non-zero distances are distinct (except, of course, for the forced equalities of the form $d(x,y)=d(y,x)$). The set of distance separated spaces of $\FM_m$ will be denoted $\ds_m$. 
\end{defn}

The following lemma is obvious.

\begin{lem}\label{l:d-sep}
The following conditions are equivalent:
\begin{enumerate}
  \item $M$ is distance separated.
  \item if $d(x,y)=d(x',y')\neq 0$, then $(x,y)=(x',y')$, or $(x,y)=(y',x')$.
\end{enumerate}
\end{lem}

In the rest of the section we show that $\ds_m$ is open and dense in $\FM_m$. Let $B_{\FM_m}(M,\eps):=\{M'\in \FM_m| d_B(M,M')<\eps\}$, denote the open ball of radius $\eps$, centered at $M$.

\begin{thm} \label{t:dsopen}
 $\ds_m$ is open in $\FM_m$. Explicitly, suppose $M\in \ds_m$ and $0<\eps<\de/10$, where $\de:=\de(DM)$. Then $B_{\FM_m}(M,\eps)\subset \ds_m$.
\end{thm}

\begin{proof}
Let $M'\in B_{\FM_m}(M,\eps)$ and choose a bijection $f:M\to M'$ with $\textrm{sep }f<\eps$. By the triangle inequality,
\[
|d(x,y)-d(f(x),f(y))|\leq 2\eps, \quad \textrm{for all } x,y\in M.
\]
Let  $x',y',u',v'\in M'$, and let $g$ denote $f^{-1}$. Then
\begin{align*}
|d(g(x'),g(y'))-d(g(u'),g(v'))|   &\leq       |d(g(x'),g(y'))-d(x',y')|+   \\
      & \quad +|d(x',y')-d(u',v')|+ \\
      & \quad\quad +|d(u',v')-d(g(u'),g(v')) \\
       &\leq  |d(x',y')-d(u',v')| +4\eps
\end{align*}
To see that $M'$ is distance separated,  suppose now that $d(x',y')=d(u',v')\neq 0$. Replacing in the above inequality, we obtain:
\[
|d(g(x'),g(y'))-d(g(u'),g(v'))| \leq 4\eps < \de.
\]
It follows that $d(g(x'),g(y'))=d(g(u'),g(v'))\neq 0$. Since $M\in \ds_m$, Lemma \ref{l:d-sep}(ii) implies, 
\[
(g(x'),g(y')) =
\begin{cases}  
     &(g(u'),g(v'))\\
     & \quad\quad \text{or} \\
      & (g(v'),g(u')) .
\end{cases}
\]
Since $g$ is a bijection, $(x',y')=(u',v')$ or $(x',y')=(v',u')$. Appealing again to Lemma \ref{l:d-sep}(ii), we see that $M'$ is distance separated, as desired.
\end{proof}

The next result shows that $\ds_m\subset \FM_m$ is dense.

\begin{thm} \label{t:dsdense}
 The set of distance separated spaces is dense in $\FM_m$. Explicitly, let $M \in\FM_m$ and $\eps>0$, be given. Then there is $M' \in\ds_m$ such that $d_B(M,M')<\eps$.
\end{thm}

\begin{proof}
By induction on $m=|M|$. The result is obvious for $m=2$. Suppose the proposition is true for every space with $m\geq 2$. Let $M=\{x_1,\dots,x_{m+1}\}\subset \R^N$ be a space with $m+1$ points, and $\de=\de(DM)$. Clearly, we can assume that $\eps <\de/10$.

Let $M_m:=\{x_1,\dots,x_m\}$ and set $\de_m:=\de(DM_m)$. Then $\eps<\de_m/10$, since $\de\leq \de_m$. By the induction hypothesis, there is a $d$-separated $M'_m\subset \R^N$ and a bijection $f_m:M_m \to M'_m$,  such that $d(x_i,x'_i)<\eps$, for $i=1,\dots,m$, where $x'_i:=f(x_i)$. Also, since $M'_m$ is $d$-separated, we have $|(DM'_m)_+| = (m-1)m/2$.

Let $\si_{kij}:=\SI(x'_k,d(x'_i,x'_j))\}$, for $k=1,\dots,m$ and $1\leq i < j \leq m$, and let $\ell_{ij}$ denote the hyperplane in $\R^N$ through the midpoint of the segment with endpoints $x'_i,x'_j$, that is perpendicular to it. Finally, define:
\[
L:= \bigcup_{k=1}^m\quad \bigcup_{1\leq i < j\leq m} (\si_{kij}\cup \ell_{ij}).
\]
Then $L\subset \R^N$ is closed and has no interior points. Hence, for any $x\in \R^N$ and $ r>0$, $B(x,r)\setminus L\neq \emptyset$.

Then, for any $x' \in B(x_{m+1}, \eps)\setminus L$ and $1\leq i,j,k \leq m,\, (i<j)$, we have:
\begin{enumerate}
  \item $d(x',x'_i)=d(x',x'_j)$ if and only if $x'\in \ell_{ij}$.
  \item if $d(x',x'_k)=d(x'_i,x'_j) $ if and only if $x'\in \si_{kij}$.
\end{enumerate}
To see (i), note that the points $x',x'_i,x'_j$ determine a plane (if they are not collinear) or a line in $\R^N$. Intersecting this plane or line with $\ell_{ij}$ reduces the claim to an obvious assertion about isosceles triangles in a plane, or about a point ($x'$) in a segment (the one determined by $x'_i,x'_j$).

For (ii), note that $d(x',x'_k)=d(x'_i,x'_j):=r'_{ij}\in DM'_+ $ means that $x'\in \SI(x_k',r'_{ij})$.

Next, we choose any $x'_{m+1} \in B(x_{m+1}, \eps)\setminus L$, and define $f:M\to M':=M'_m\cup \{x'_{m+1}\}$ by $f|_{M_m}:= f_m$ and $f(x_{m+1}):=x'_{m+1}$. Then $M'\in \FM_{m+1}$, $f$ is a bijection and, clearly, $d_B(M,M')<\eps$. It remains to show that $M'$ is $d$-separated. Given that $|(DM'_m)_+|=(m-1)m/2$, to prove that $|(DM'_{m+1})_+|=m(m+1)/2$, it suffices to see that the $m$ distances $d(x'_{m+1},x'_k), (1\leq k \leq m)$, are $\neq 0$ (true by the choice of $x'_{m+1}$ and $\eps$), are different among themselves (true by (i) above) and that they are different from all of the $d(x'_i,x'_j)\in (DM_m)_+$ (true by (ii)). This concludes the proof.
\end{proof}

\section{The ``general case'' for $CS(M)$ }
\label{s:gen_case}
This section concerns stability under small perturbations of the initial data. Example  \ref{exa:unst} shows that neither $CS(M)$ nor $MC(M)$ are stable constructions nor need they be trees. We show, however, that producing a small perturbation of $M$, if necessary, will guarantee that $CS(M)$ is a tree, a result with important consequences in practice because, in many cases, it guarantees effective computation (for instance, of the finite dimension of $M$, via $CS(M)$, as explained in Section \ref{s:mtogr}), since a slight modification of $M$ is within empirical error. No such result holds for $MC(M)$.

Theorem \ref{t:CStree} says that the ``general'' case is for $CS(M)$ to be a tree, and justifies the name: while $CS(M)$ need not be sparse, we can make it sparsest possible, a tree, at the price of modifying $M$ slightly, as little as we wish.

Recall the notation and construction of Section \ref{s:csg}. We assume that $|K_i|>1$, that no $k\in K_i$ has cycles, and that $M$ is distance separated. We take a closer look at cycles in $S_{i+1}$ for $i\geq 0$.

\begin{figure}
\begin{tikzpicture}
\draw[rotate= 40] (0,0) ellipse (30pt and 15pt);
\draw[rotate= 0] (2,1.5) ellipse (30pt and 15pt);
\draw[rotate= 0] (5,1.5) ellipse (30pt and 15pt);
\draw[rotate= 0] (2,-1.5) ellipse (30pt and 15pt);
\draw[rotate= 0] (5,-1.5) ellipse (30pt and 15pt);
\draw[dashed] (7,1.5) -- (9,1.1);
\draw[dashed] (7,-1.5) -- (9,-1.1);
\filldraw[black] (-.6,-.6) circle (1pt) node[anchor=south]{$x_1$};
\filldraw[black] (-.1,-.4) circle (1pt) node[anchor=west]{$x_2$};
\draw (-.6,.8) node[anchor=north] {$k_1$};
\draw[gray] (-.6,-.6) -- (-.1,-.4);
\filldraw[black] (-.4,-.1) circle (1pt) ;
\draw[gray] (-.1,-.4) -- (-.4,-.1);
\draw[dashed, gray]  (-.4,-.1) -- (0,.3);
\draw[dashed, gray]  (.4,-.1) -- (0.6,.6);
\draw[dashed, gray]  (0,.3) -- (.4,-.1);
\filldraw[black] (.6,.6) circle (1pt) node[anchor=west]{$x_{(j_1-1)}$};
\draw (1.3,2.5) node[anchor=north] {$k_2$};
\draw (.95,1.2) node[anchor=east] {$e_1$};
\filldraw[black] (1.05,1.5) circle (1pt) node[anchor=west]{$x_{j_1}$};
\draw[very thick]  (0.6,.6) -- (1.05,1.5);
\draw[dashed, gray]  (1.05,1.5) -- (1.9,1.1);
\draw[dashed, gray]  (1.9,1.1) -- (2.1,1.9);
\filldraw[black] (2.9,1.5) circle (1pt) node[anchor=north]{$x_{(j_2-1)}$};
\draw[dashed, gray]  (2.1,1.9) -- (2.9,1.5);
\draw[very thick]  (2.9,1.5) -- (4.3,1.8);
\filldraw[black] (2.9,1.5) circle (1pt) node[anchor=north]{$x_{(j_2-1)}$};
\draw (3.85,1.9) node[anchor=east] {$e_2$};
\filldraw[black] (4.3,1.8) circle (1pt) node[anchor=north]{$x_{j_2}$};
\draw (6,2.3) node[anchor=north] {$k_3$};
\draw[dashed, gray] (4.3,1.8) -- (5.5,1.4);
\draw[very thick] (-.6,-.6) -- (1.2,-1.5);
\filldraw[black] (1.2,-1.5) circle (1pt) node[anchor=west]{$x_n$};
\draw (1.9,-.4) node[anchor=north] {$k_p$};
\draw (4,-1.1) node[anchor=east] {$e_{p-1}$};
\draw[very thick]  (2.7,-1.35) -- (4.6,-1.2);
\filldraw[black] (2.7,-1.35) circle (1pt) node[anchor=north]{$x_{j_p}$};
\filldraw[black] (4.6,-1.2) circle (1pt) node[anchor=north]{$x_{(j_p-1)}$};\draw[dashed, gray]  (4.6,-1.2) -- (5.5,-1.3);
\draw[dashed, gray]  (1.8,-1.2) -- (2.7,-1.35);
\draw[dashed, gray]  (2.1,-1.9) -- (1.8,-1.2);
\draw[dashed, gray]  (2.1,-1.9) -- (1.2,-1.5);
\draw (5.9,-.4) node[anchor=north] {$k_{p-1}$};
\draw (.2,-1.1) node[anchor=north] {$e_p$};
\draw[dashed, gray]  (9.5,.7) -- (10,0) -- (9.5,-.7);
\end{tikzpicture}
\caption{Illustration of Lemma \ref{l:i-edges}:  a typical cycle $z$ of  $S_{i+1}$, where $x_1,\dots,x_n$ are the vertices of $z$, $k_1,\dots,k_p$ are the components of $S_i$, and   $e_1,\dots,e_p$ are the ($i+1$)-edges.}
\label{fig:cycle}
\end{figure}
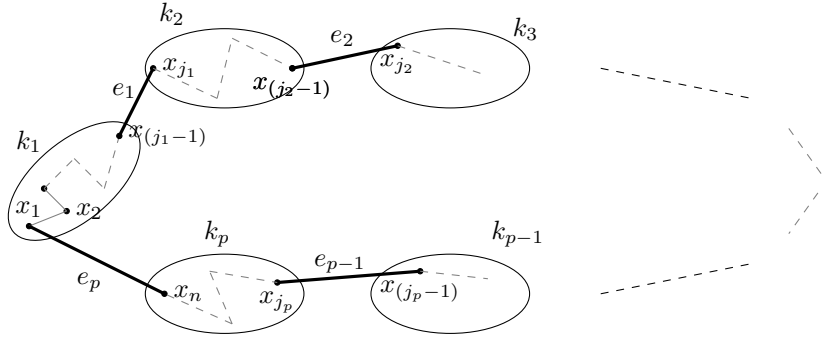

Let $z=x_1,\dots,x_n$ be a cycle in $S_{i+1}$. We wish to express $z$ in terms of the components of $S_{i}$ (cf. Fig. \ref{fig:cycle}). Cyclically renumbering the vertices of $z$, if necessary, we may assume that $[x_1]_i\neq [x_n]_i$. Suppose
\begin{align}
\label{}
    k_1:=[x_1]_i=& \{x_1,\dots, x_{(j_1-1)}\} \nonumber\\
   k_2:=[x_{j_1}]_i=&  \{x_{j_1},\dots,x_{(j_2-1)} \} \\
   \vdots    \nonumber \\
   k_p:=[x_{j_{p}}]_i=& \{x_{j_{p}},\dots,x_n \} \nonumber
\end{align}
Consecutive components $k_h,k_{h+1}$ are connected (in $S_{i+1}$) by edges:
\begin{equation}
e_h =  \{x_{(j_h-1)} ,x_{j_h}\},\, (1\leq h\leq p-1), \quad e_p=\{x_n,x_1\}
\end{equation}
that we call \emph{$(i+1)$-edges}, for short.

We summarize our results in a lemma:

\begin{lem}\label{l:i-edges}
If $M\in \ds_m$ and no component of $S_i$ has cycles, then any cycle $z$ of $S_{i+1}$ can be written in the form
\[
z=k_1e_1k_2\dots e_{p-1}k_pe_p
\]
where $k_j\in K_i$, $e_j\in F(S_{i+1})$, $(j=1,\dots, p)$ and $p\geq 3$. Moreover, two different components can be connected by at most one $(i+1)$-edge.
\end{lem}
\begin{proof}
Suppose $k, k'\in K_i$ are different, $x,y\in k$, $x',y'\in k'$, and both $\{x, x'\},\{y, y'\}$ belong to $FS_{i+1}$. By Prop. \ref{p:tech}, $d(x, x')= \max\{\nu_i(k),\nu_i(k')\}=d(y, y')$. By Lemma \ref{l:d-sep}, $(x, x')=(y, y')$ or  $(x, x')=(y', y)$. The first possibility leads to the claim, as desired, while the second leads to $k=k'$, a contradiction.

The hypothesis that no component of $S_i$ has cycles implies that $z$ cannot be contained in a single component, thus $p\geq 2$. If $p=2$ i.e., if $z=k_1e_1k_2e_2$, then $e_1=e_2$. So $\{x_{(j_1-1)} , x_{j_1} \}=\{x_n , x_1\}$, and we consider two cases. If $x_{(j_1-1)} =x_n$, then $k=k'$, a contradiction. The other possibility implies $k=\{x_1\}$ and $k'=\{x_2\}$, so that $n=2$ and $z$ is not a cycle. Hence $p\geq 3$, as desired.
\end{proof}

\begin{thm} \label{t:tree}
If $M$ is $d$-separated, then $CS(M)$ is a tree.
\end{thm}

\begin{proof}
We prove, by induction, that the components of $S_i$ have no cycles. This is clear for $i=0$, since the components are single points. Assume the result is true for some $i\geq 0$. If $|K_i|=1$, then $i\geq 1$ and $S_i=CS(M)$ is a tree. If $|K_i|>1$, then we apply step $i+1$ to construct $S_{i+1}$. 

Suppose, for contradiction, that $S_{i+1}$ has a cycle $z=x_1,\dots,x_n=k_1e_1\dots,k_pe_p$, where $p\geq 3$, $k_h\in K_i$ and $e_h\in FS_{i+1}$, according to Lemma \ref{l:i-edges}. Since $M$ is distance separated, the lengths $\ell(e_h)$ are distinct. Suppose, without loss of generality, that $e_2$ is the longest $(i+1)$-edge of $z$. In particular, $\ell(e_2)>\max\{\ell(e_1),\ell(e_3)\}$. By Prop. \ref{p:tech}, $\ell(e_2)=\max \{ \nu_i(k_2),\nu_i(k_3)\}$. On the other hand, $\ell(e_1)\geq \nu_i(k_2)$ and $\ell(e_3)\geq \nu_i(k_3)$. Thus,
\[
\ell(e_2)=\max \{ \nu_i(k_2),\nu_i(k_3)\} \leq \max \{ \ell(e_1),\ell(e_3)\} < \ell(e_2),
\]
a contradiction. The proof is complete.
\end{proof}

Putting together theorems \ref{t:dsopen}, \ref{t:dsdense} and \ref{t:tree}, we obtain the main result of the paper:

\begin{thm} \label{t:CStree}
Let $M \in \FM_m$, $ \de=\de(DM)$, and $\eps<\de/10$. Then there is $M' \in \ds_m$, such that $d_B(M, M')<\eps$ and $CS(M')$ is a tree. If, moreover, $M\in \ds_m$, then $CS(M')$ is a tree, for all $M'\in B_{\FM_m}(M,\eps)$.
\end{thm}

%\newpage
%XXXXXXXXXXXXXXXXXXXXXX
%XXXXXXXXXXXXXXXXXXXXXX
%XXXXXXXXXXXXXXXXXXXXXX

\section{Invariance}
\label{s:inv}
In this section we prove that the constructions are invariant under isometries.

\begin{thm} \label{t:isometry}
Let $(M,d),(M', d')$ be finite metric spaces. Then any isometry $\varphi:(M,d)\to (M',d')$ induces graph isomorphisms $MC(M) \cong MC(M')$, and $CS(M) \cong CS(M')$.
\end{thm}

\begin{proof}
We first show that $\varphi$ defines a graph isomorphism $\varphi: MC(M) \cong MC(M')$. Since $d(x,y)\leq r$ if and only if $d'(\varphi(x),\varphi(y)) \leq r$, we have $\{x,y\}\in G_r$ if and only if $\{\varphi(x),\varphi(y)\} \in G'_r := \{(x',y')\in M'\times M' | d'(x',y')\leq r  \}$. Thus $\varphi: MC(M) \cong MC(M')$.

Next, consider the construction of $CS(M), CS(M')$. Consider the diagram: 
\[
\begin{tikzcd}[]
S_0 \arrow[d, "\varphi"] \arrow[r, hook,"\textrm{inc}"]  & S_1 \arrow[d, "\varphi"]  \arrow[r, hook,"\textrm{inc}"] &\cdots \arrow[r, hook,"\textrm{inc}"]  &  S_i   \arrow[d, "\varphi"]  \arrow[r, hook,"\textrm{inc}"]  & S_{i +1}  \arrow[d, dashrightarrow, "\varphi"]\\
S'_0  \arrow[r, hook,"\textrm{inc}"]      & S'_1   \arrow[r, hook,"\textrm{inc}"] & \cdots \arrow[r, hook,"\textrm{inc}"] & S'_i \arrow[r, hook,"\textrm{inc}"] & S'_{i +1}
\end{tikzcd}
\]
where the vertical arrows restricted to $V(S_i)=M\to V(S'_i)=M',\,i \geq 0$, are all equal to $\varphi$, so that the diagram commutes trivially. We assume $\varphi:S_j \to S'_j$, for $0\leq j \leq i$, is an isomorphism, and prove that also $\varphi:S_{i+1} \to S'_{i+1}$ is an isomorphism, starting with the case $i=0$. 

We claim that $\nu_0:M\to \R_+$, equals $\nu'_0\varphi$, where $\nu'_0:M'\to \R_+$ is the function on $M'$ defined analogously to $\nu_0$. To see this, note that $\nu_0(x)=d(x,y_0)$, for some $y_0\in M$. Then $d'(\varphi(x), \varphi(y_0)) =   d(x,y_0)$, so that $\nu'_0(\varphi(x)) \leq \nu_0(x)$. The reverse inequality is proved similarly. 

To see that $\varphi:S_1\to S'_1$ is an isomorphism of graphs, we need to check that $x$ is adjacent to $y$ iff $\varphi(x)$ is adjacent to $ \varphi(y)$, for all $x,y\in M=V(S_1)$, By definition, $x$ is adjacent to $y$ iff $d(x,y)=\nu_0(x) = \nu'_0(\varphi(x)) =d'(\varphi(x),\varphi(y))$, which is precisely the condition for $\varphi(x)$ to be adjacent to $\varphi(y)$. 

The proof that $\varphi:S_{i + 1}\to S'_{i + 1}$, is an isomorphism (provided that $\varphi:S_{i }\to S'_{i}$ is) is entirely similar, once we have checked that $\nu_i= \nu'_i \varphi:K_i \to \R_+$. If $x\in k\in K_i$, we claim that $\varphi(k)$ is the connected component of $\varphi(x)$ in $S'_i$. This is clear, given the induction hypothesis: any edge $e\in E(S_i)$ is mapped to an edge $\varphi(e)\in E(S'_i)$. Hence, any path joining $x$ to some $y$, will be mapped by $\varphi$ to a path in $S'_i$, and similarly for $\varphi^{-1}$. Thus, $\varphi$ is a bijection from the path components of $S_i$ to the path components of $S'_i$, and is an isomorphism on each of them. For $k\in K_i$, let $\nu_i(k)=d(x_0,y_0)$, for certain $x_0\in k$, and $y_0\in \overline{k}$. Then $d(x_0,y_0)=d'(\varphi(x_0),\varphi(y_0))$, where $\varphi(y_0)\notin\varphi(k)$, for otherwise $y_0\in k$. It follows that $\nu'_i(\varphi(k)) \leq \nu_i(k)$, as desired. The reverse inequality is proved similarly. 

Again, to see that $\varphi:S_{i +1}\to S'_{i +1}$ is an isomorphism, it suffices to show that $x$ and $ y$ are adjacent iff $\varphi(x)$ and $\varphi(y)$ are adjacent, for $\{x,y\}\in F(S_{i+1})$. By definition, this means that $d(x,y)=\nu_i(x)$, for some $x\in k,\, y\notin k$. Then $d'(\varphi(x),\varphi(y))=d(x,y)=\nu_i(x)=\nu'_i(\varphi(x))$, as desired. This completes the proof of the theorem.
\end{proof}

\section{Length spaces and graphs}
\label{s:intri}

What should it mean to say that $(M,d)$ is a \emph{length space} or, equivalently, that $d$ is a \emph{path-metric}? Our definition is, roughly speaking, that $M$ is a length space if there is a graph (with $M$ as vertex set and $w=d$) whose path-metric coincides with $d$. Here are the details. We start by considering the set of all (connected) graphs whose path-metric is $d$, $\II:=\{G^w| V(G)=M \wedge d_G^w=d\}$, where $d_G^w$ is defined using the weights $w_{\{ x, y \}} := d(x,y)$. As metric spaces, $(M,d)=(G^w,d_G^w)$, for all $G^w\in\II$ (more precisely, $(M,d)=(V(G),d_G^w)$, but we will continue to abuse language freely). Note that $\II\neq \emptyset$, since $K_{|M|}^w\in \II$. To simplify the notation, we will not always decorate graphs of $\II$ with the weights $w$.

Consider the subset $Q\subset E(K_{|M|})$, defined by: 
\begin{displaymath}
Q:=\Big\{ \{x,y\}\in E(K_{|M|})|   \textrm{ there is a geodesic } p(x,y) \textrm{ of $K_{|M|}$, with } c(p)\ge 2   \Big\}.
\end{displaymath}

\begin{prop}
Considered as a poset with respect to inclusion, $\II$ has a maximum and a minimum. The maximum is $K_{|M|}^w$, and the minimum, denoted $\SI_M^w$, is the graph with edges  $E(\SI_M):=E(K_{|M|})\setminus Q$.
\end{prop}

\begin{proof} Note that the maximum is the union of all $G\in \II$, and the minimum is the intersection of all $G\in\II$, provided it belongs to $\II$. Also, since all graphs $G \in \II$ have vertex set $M$, they differ only on their sets of edges. Thus, a union or intersection of graphs in $\II$ is just the union or intersection of their edges. 

The maximum is obviously equal to $K_{|M|}$. To identify the minimum, it suffices to show that $\SI_M$ belongs to $\II$, and is contained in the intersection of all graphs of $\II$. Let's show first that $\SI_M\in \II$. Given a pair $x,y\in M$, among all geodesics joining $x,y$ in $K_{|M|}$, let $p(x,y)=e_1\dots e_n$, be one with \emph{largest} count. If all $e_i\in \SI_M$, then $\SI\in \II$, since $d_{\SI_M}(x,y)=\ell(p)=d(x,y)$. If, on the contrary, some $e_i\in Q$, we can replace $e_i$ by a geodesic (in $K_{|M|}$) $p'$, with $c(p')>1$. Then the new geodesic $e_1\dots e_{i-1}p' e_{i+1}\dots e_n$ has count $> n+1$, a contradiction. Thus, $\SI_M\in \II$.

Next, we show that if $e=\{ x, y \}\notin Q$, then $e$ belongs to every graph $G\in \II$. Indeed, all geodesics $p(x,y)$ must have count 1, i.e. $p=e$. Hence $\{ x, y \}\in E(G)$, as desired.
\end{proof}

Note that $\SI_M$ is particularly interesting because it is the smallest graph on $M$ whose path-metric is $d$. To this feature of obvious computational importance, we add another, more theoretical one:

\begin{defn}
We say that $(M,d)$ is an \emph{intrinsic} finite metric space, or that $d$ is \emph{intrinsic} (equivalently, $M$ is a \emph{path-space}, or $d$ is a \emph{path-metric}) if $\SI_M\neq K_{|M|}$; equivalently, if $Q\neq \emptyset$. A metric that is not intrinsic is called \emph{extrinsic}; equivalently, if $\SI_M = K_{|M|}$ which is equivalent, in its turn, to $Q= \emptyset$.
\end{defn}

By way of motivation, observe that every $(M,d)$ has a graph whose path-metric is $d$, namely $K_{|M|}$. So we require that $\SI_M$ be strictly smaller than $K_{M}$. This condition guarantees that at least some of the geodesics defining the path-metric will be ``honest'' geodesics, in the sense that they will have count $\geq 2$.

\begin{defn}
Suppose $(M,d)$ is intrinsic. We call $M$ \emph{intrinsic-I} if there is a constant $r>0$ such that all edges of $\Sigma_M$ have length $r$. It is called \emph{intrinsic-II} if it is not intrinsic-I. Equivalently, it is intrinsic-II if the edges of $\Sigma_M$ have non-constant length.
\end{defn}

\begin{exa}
(a) A simple example of extrinsic space is given by the points $M=\{a,b,c\}$ of the plane, with coordinates $(0,0),\, (1,0), \, (1,1)$, respectively. The metric is either $\ell_1,\ell_2$ or $\ell_\infty$, so that $d(a,b)=d(b,c)=1$, and $d(a,c)=2,\sqrt{2}$, or $ 1$, respectively. Then $M$ is a path-space when the metric is induced by $\ell_1$. In the other two cases $Q=\emptyset$, and the metric is extrinsic. With $\ell_1$, $M$ is intrinsic-I, since the two edges of $\SI_M^w$ have $r=1$.

(b) This is an example of an intrinsic-II space. Let $M=\{x,y,z,t\}\subset\R^2$, where $x=(0,0), \, y=(1,0), \, z=(2,0), \, t=(1,1)$, with the induced euclidean distance. Then $\Sigma_M^w$ has edges (with weights) $xy=\{x,y,1\}$, $xt=\{x,t,\sqrt{2}\}$, $yz=\{y, z,1\}$, $yt=\{y, t,1\}$, and $tz=\{t, z,\sqrt{2}\}$. Note that $CS(M)=MC(M)$ has edges $xy,yz$ and $yt$.
\end{exa}

\section{Relations among the graphs}

\begin{figure}[h!]
\begin{tikzpicture}
  \node (max) at (0,3) {$K_{|M|}$};
  \node (a) at (-2,0) {$\Sigma_M$};
  \node (b) at (2,1.5) {$MC(M)$};
  \node (c) at (0,-1.5) {$\Sigma_M \cap MC(M)$};
  \node (min) at (0,-3) {$CS(M)$};
  \draw (min) -- (c) -- (a) -- (max)
  (c) -- (b) -- (max);
\end{tikzpicture}
\caption{Hasse diagram of inclusions}
\label{fig:hasse}
\end{figure}
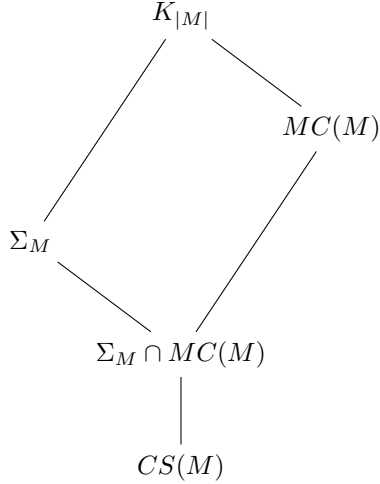

In this section we make explicit several relations that hold among the various graphs defined in the paper, summarized in Fig. \ref{fig:hasse}.

\begin{thm} \label{t:rela2}
The relationship among the various graphs considered here is summarized in the following assertions:
\begin{enumerate}
  \item The graph $CS(M)$ is a subgraph of both $\SI_M$ and $MC(M)$.
  \item If $M$ is intrinsic-I, then $\SI_M =CS(M) = MC(M)$.
  \item If $M$ is intrinsic-II and $\Sigma_M$ is a tree, then $\SI_M =CS(M)$. In general, $\Sigma_M$, $CS(M)$ and $MC(M)$ are distinct (cf. Example \ref{exa:exa2} below).
\end{enumerate}
\end{thm}

\begin{proof}
(i) We begin by proving inclusion in $\SI_M$. Let $e=\{x,y\}\in E(S_1)$, and suppose, for contradiction,  that $\{x, y\}\in Q$. By definition, we can find $z\in M$, such that $d(x,y)=d(x,z)+d(z,y)$. Thus $d(x,z)<d(x,y)=\nu_0(x)$, a contradiction, hence   $S_1\subseteq \SI_M$. A similar argument applies to the inductive step: assuming  $S_i\subseteq \SI_M$, if for some $e=\{x,y\} \in F(S_{i +1}),\,(x\in k)$, we have $e\in Q$, then, as before, there is $z$ with $d(x,z)<d(x,y)=\nu_i(k)$. Hence $z\in k$, and $d(z,y)<\nu_i(k)$, a contradiction, since $\nu_i(k)$ is the minimum distance between $k$ and the connected component of $y$.

For the second inclusion, suppose $MC(M)=G_t$ is defined by a cut value $t\geq \nabla$, where $\nabla:= \max\{\nu_0(x)|x \in M\}$. If $e=\{x,y\} \in E(S_1)$, then $d(x,y)=\nu_0(x)\leq \nabla$, so that $e\in G_t$. Suppose now that $S_i\subseteq G_t$, and let $e=\{x,y\}\in F(S_{i +1})$, where $x\in k$ and $y\in k_1\neq k$. Since $k,k_1$ are both contained in the connected $G_t$, there is a path, which we can, and do, assume there is a geodesic, $q(x,y)=e_1\dots e_n$ in $G_t$, with $n\geq 1$, and, say,  $e_1=\{x,x_1\}$. We consider first the case where $x_1\notin k$. Then, $d(x,y)=\nu_i(x)\leq d(x,x_1)$ and, since $e_1\in G_t$, $d(x,x_1)\leq t$. It follows that $e\in G_t$, as desired. We leave it to the reader to reduce to this case, if $x_1\in k$. This completes the proof of (i).

(ii) The hypothesis mean that $(M,d)=(\Sigma_M,d_\Sigma)$, where the edges of $\Sigma_M$ have constant length $r>0$. From the definition of $\Sigma_M$ it follows that $e=\{x,y\}$ belongs to $\Sigma_M$ iff $d(x,y)=r$. This implies immediately that $\Sigma_M=S_1= CS(M)$. It is also clear that any cut-value $t<r$ will result in a graph with no edges, while $t=r$ will give the connected graph $\Sigma_M$ which, hence, equals $MC(M)$, as desired.

%\vspace{2mm}
 (iii) By (i), $CS(M)\subseteq \Sigma_M $. Hence, $CS(M)$ is a tree and, since both graphs have the same set of vertices, $CS(M)$ must necessarily coincide with $\Sigma_M$. 
\end{proof}

It follows from the theorem that $CS(M)$ is obtained from $\SI_M$ by removing edges and is, in general, strictly smaller than $\SI_M$. We thus have:

\begin{cor}
Either $CS(M)$ distorts the metric $d$ (when $CS(M)\subsetneq \SI_M$), or it is the smallest graph (with vertex set $M$) that gives $M$ the structure of a length space (when $CS(M)= \SI_M$).
\end{cor}

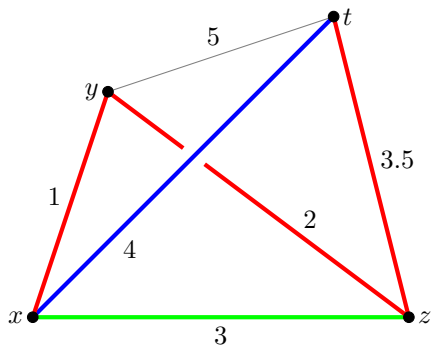
\begin{figure}
\begin{tikzpicture}[] 
\draw[blue, ultra thick] (4,4) -- (0,0);
\draw[green, ultra thick] (0,0)  -- (5,0);
\draw[red, ultra thick] (0,0)  -- (1,3)  (5,0) -- (4,4);
\draw[red, ultra thick] (1,3) -- (2,9/4)   (16/7,57/28) -- (5,0);
\draw[gray] (1,3)  -- (4,4);
\filldraw[black] (0,0) circle (2pt) node[anchor=east]{$x$};
\filldraw[black] (1,3) circle (2pt) node[anchor=east]{$y$};
\filldraw[black] (5,0) circle (2pt) node[anchor=west]{$z$};
\filldraw[black] (4,4) circle (2pt) node[anchor=west]{$t$};
\draw (2.5,0) node[anchor=north] {$3$};
\draw (0.5,1.6) node[anchor=east] {$1$};
\draw (2.4,3.5) node[anchor=south] {$5$};
\draw (4.5,2.1) node[anchor=west] {$3.5$};
\draw (3.9,1.3) node[anchor=east] {$2$};
\draw (1.5,.9) node[anchor=east] {$4$};
\end{tikzpicture}
\caption{$K_{|M|}$ with weights. $CS$ in red, $\SI_M$ equals $CS$ plus the blue edge ($xt$) and $MC$ equals $CS$ plus the green edge ($xz$)}
\label{fig:exa2}
\end{figure}

\begin{exa} \label{exa:exa2}
Consider the intrinsic-II space $M=\{x,y,z,t\}$ (see Fig. \ref{fig:exa2}), with $E(K_{|M|})=\{ xy=(x,y,1), xz=(x,z,3), xt=(x,t,4), yz=(y, z,2), yt=(y, t,5), zt=(z,t,3)    \}$. Then $E(\Sigma_M)=\{xy, xt, yz, zt \}$, $E(MC(M))=\{xy, xz, yz, zt\}$, and $E(CS(M))=\{xy, yz, zt\}$. In this example, $CS(M)= \Sigma_M \cap MC(M)  \subsetneq \Sigma_M$ and $CS(M)\subsetneq MC(M)$. Both $CS(M)$ and $MC(M)$ distort distance. 
\end{exa}

%\section{Bibliography}

%\bibliographystyle{amsplain}
%\bibliography{xbib}

%\bibliography{}

% -----------------------------------------------------------

% -----------------------------------------------------------
%\bibliographystyle{amsplain}
%\bibliography{xbib}
%FinDIMplants.bib

\end{document}